\documentclass[12pt]{article}

\begin{document}

\begin{center}
{\Large
\textbf{Thoughts on the Riemann Hypothesis}
}
\vspace{6mm}
\\
G. J. Chaitin, \emph{IBM Research}
\\
{\footnotesize
\texttt{http://www.cs.auckland.ac.nz/CDMTCS/chaitin}
}
\end{center}

\section*{}

The simultaneous appearance in May 2003 of \textbf{four} books
on the Riemann hypothesis (RH) provoked these reflections.
We briefly discuss whether the RH should be added as a new axiom, or
whether a proof of the RH might involve the notion of randomness.

\section*{New pragmatically-justified mathematical axioms that are not at all self-evident}

A pragmatically-justified principle is one that is justified by its many important 
\textbf{consequences}---which is precisely the opposite of normal mathematical practice.\footnote
{However new mathematical \textbf{concepts} such as $\sqrt{-1}$ and Turing's definition of
computability certainly are judged by their fruitfulness\emph{---Fran\c{c}oise 
Chaitin-Chatelin, private communication.}}
However this is standard operating procedure in physics.

Are there mathematical propositions for which there is a
considerable amount of computational evidence, evidence that is so persuasive that a physicist
would regard them as experimentally verified?
And are these propositions fruitful? Do they yield many other signficant results?

Yes, I think so.
Currently, the two best candidates\footnote
{Yet another class of pragmatically-justified axioms are the large cardinal axioms 
or the axiom of determinacy used in set theory, as discussed in Mary Tiles,
\emph{The Philosophy of Set Theory,} Chapters 8 and 9.
For the latest developments, see Hugh Woodin, ``The continuum hypothesis,''
\emph{AMS Notices} \textbf{48} (2001), pp.\ 567--576, 681--690.}  
for useful new axioms of the
kind that G\"odel and I propose [1]
that are justified pragmatically as in physics are: 
\begin{itemize}
\item
the $\mathbf{P} \neq \mathbf{NP}$ hypothesis in theoretical computer science
that conjectures that many problems require an exponential amount of work to resolve,
and 
\item
the
Riemann hypothesis concerning the location of the complex
zeroes of the Riemann zeta-function
\[
   \zeta(s) = \sum_n \frac{1}{n^s} = \prod_p \frac{1}{1-\frac{1}{p^s}} .
\]
(Here $n$ ranges over positive integers and $p$ ranges over the primes.)\footnote
{You start with this formula and then you get the full zeta function by analytic continuation.}
\end{itemize}
Knowing the zeroes of the zeta function, 
i.e., the values of $s$ for which $\zeta(s) = 0$,
tells us a lot about the smoothness of the distribution of prime numbers, 
as is explained in these four books:
\begin{itemize}
\item
Marcus du Sautoy, 
\emph{The Music of the Primes,} Harper Collins, 2003.
\item
John Derbyshire, 
\emph{Prime Obsession,} Joseph Henry Press, 2003.
\item
Karl Sabbagh, 
\emph{The Riemann Hypothesis,} Farrar, Strauss and Giroux, 2003.
\item
Julian Havil, 
\emph{Gamma,} Princeton University Press, 2003.\footnote
{Supposedly Havil's book is on Euler's constant $\gamma$, not the RH, but ignore that.
Sections 15.6, 16.8 and 16.13 of his book are particularly relevant to this paper.}
\end{itemize}

The Riemann zeta function is like my $\Omega$ number: it captures a lot of information
about the primes in one tidy package. $\Omega$ is a single real number that 
contains a lot of information about the halting problem.\footnote
{$\Omega = \sum_{\mbox{\tiny $p$ halts}} 2^{-|p|}$ 
is the halting probability of a suitably chosen universal Turing machine.
$\Omega$ is ``incompressible'' or ``algorithmically random.'' 
Given the first $N$ bits of the base-two expansion of $\Omega$, one can determine whether each
binary program $p$ of size $|p| \le N$ halts.
This information cannot be packaged more concisely.
See [2], Sections 2.5 through 2.11.}
And the RH is useful because it contains a lot of number-theoretic information:
many number-theoretic results follow from it.
 
Of the authors of the above four books on the RH, the one who takes G\"odel most
seriously is du Sautoy, who has an entire chapter on G\"odel and Turing in his book.  
In that chapter on p.\ 181,
du Sautoy raises the issue of whether the RH might require new axioms. 
On p.\ 182 he quotes G\"odel,\footnote
{Unfortunately du Sautoy does not identify the source of his G\"odel quote.
I have been unable to find it in G\"odel's \emph{Collected Works.}} 
who specifically mentions 
that this might be the case for
the RH. And on p.\ 202 of that chapter du Sautoy points out that
if the RH is undecidable this implies that it's true, because if the RH were false
it would be easy to confirm that a particular zero of the zeta function is in the wrong place.

Later in his book, on pp.\ 256--257,
du Sautoy again touches upon the issue of whether the RH might require a new axiom.
He relates how Hugh Montgomery sought reassurance from G\"odel that a famous number-theoretic
conjecture---it was the twin prime conjecture, which asserts that there are infinitely many pairs 
$p$, $p + 2$
that are both prime---does not require new
axioms.
G\"odel, however, was not sure. In du Sautoy's words, sometimes one needs
``a new foundation stone to extend the base of the edifice'' of mathematics,
and this might conceivably be the case both for the twin prime conjecture and for the RH.

On the other hand, on pp.\ 128--131 du Sautoy tells the story of the Skewes number, 
an \textbf{enormous} number 
\[
   10^{10^{10^{34}}}
\]
that turned up in a proof that
an important conjecture must fail for \textbf{extremely large} cases.
The conjecture in question 
was Gauss's conjecture that the logarithmic integral 
\[
   \mbox{Li}(x) = \int_2^x \frac{du}{\ln u} 
\]
is always greater than the
number $\pi(x)$ of primes less than or equal to $x$. This was verified
by direct computation for all $x$ up to very large values. It was then refuted by
Littlewood without exhibiting a counter-example, and finally by Skewes with his enormous upper
bound on a counter-example. This raised the horrendous possibility that even though 
Gauss's conjecture is wrong, we
might \textbf{never ever} see a specific counter-example.
In other words, we might never ever know a specific value of $x$ 
for which Li($x$) is less than $\pi(x)$.
This would seem to pull the rug out from under all mathematical experimentation and 
computational evidence! However, I don't believe that it actually does.

The traditional view held by most mathematicians is that these two assertions,
$\mathbf{P} \neq \mathbf{NP}$ and the RH,
\textbf{cannot} be taken as new axioms, 
and \textbf{cannot} require new axioms,
we simply \textbf{must work much harder} to prove them. 
According to the received view, we're not clever enough, we haven't
come up with the right approach yet. 
This is very much the current consensus.
However this majority view \textbf{completely ignores}\footnote
{As du Sautoy puts it, p.\ 181, ``mathematicians consoled themselves with the belief that
anything that is really important should be provable, that it is only tortuous statements
with no valuable mathematical content that will end up being one of G\"odel's 
unprovable statements.''}
the incompleteness phenomenon discovered by G\"odel, by Turing, and extended
by my own work [2]
on information-theoretic incompleteness.
What if \textbf{there is no proof?}

In fact, new axioms \textbf{can never} be proved; if they can, they're theorems, not axioms.
So they must either be justified by direct, primordial mathematical \textbf{intuition},
or \textbf{pragmatically}, because of their rich and important consequences, as is done in physics.
And in line with du Sautoy's observation on p.\ 202, 
one cannot demand a proof that the RH is undecidable
before being willing to add it as a new axiom, because such a proof would in fact yield the
immediate corollary that the RH is true.
So proving that the RH is undecidable is no easier than proving the RH,
and the need to add the RH as a new axiom must remain a matter of faith.
The mathematical community will \textbf{never} be convinced!\footnote
{The situation with respect to 
$\mathbf{P} \neq \mathbf{NP}$ may be different. In a paper 
``Consequences of an exotic definition for $\mathbf{P} = \mathbf{NP}$,''
\emph{Applied Mathematics and Computation} \textbf{145} (2003), pp.\ 655--665, 
N. C. A. da Costa and F. A. Doria show that if ZFC 
(Zermelo-Fraenkel set theory + the axiom of choice) is consistent,
then a version of $\mathbf{P} = \mathbf{NP}$ is consistent with ZFC, so 
a version of $\mathbf{P} \neq \mathbf{NP}$ cannot be demonstrated within ZFC.
See also T. Okamoto, R. Kashima, ``Resource bounded unprovability of computational lower bounds,''
http://eprint.iacr.org/2003/187/.} 

Someone recently asked me, ``What's wrong with calling the RH a hypothesis? Why does
it have to be called an axiom? What do you gain by doing that?''
Yes, but that's beside the point, that's not the real
issue. The real question is, \textbf{Where does new mathematical knowledge come from?}

By ``new knowledge'' I mean something that cannot be deduced from our previous knowledge,
from what we already know.

As I have been insinuating, I believe that the answer to this fundamental question is that  
new mathematical knowledge comes from these three sources:
\begin{itemize}
\item[a)]
mathematical intuition and imagination ($\sqrt{-1}$!),
\item[b)]
conjectures based on computational evidence (explains calculations), and
\item[c)]
principles with pragmatic justification, i.e., rich in consequences 
(explains other theorems).\footnote
{A possible \textbf{fourth} source of mathematical knowledge is: 
d) probabilistic or statistical evidence, i.e.,
a mathematical assertion that is deemed to be true because the probability 
that it's false is immensely small, say $< 10^{-99999}$.

Here is a practical example of this: The fast primality testing algorithm currently
used in \textsl{Mathematica} does not necessarily give the correct answer, but mistakes
are highly unlikely.
Algorithms of this sort are called Monte Carlo algorithms.}
\end{itemize}
And items (b) and (c) are much like physics, if you replace ``computational evidence''
by ``experimental evidence.'' In other words, our computations are our experiments;
the empirical basis of science is in the lab, the empirical basis of math is in the computer.

Yes, I agree, mathematics and physics are different,
but perhaps they are \textbf{not} as different as most people think, perhaps it's  
a continuum of possibilities.
At one end, rigorous proofs, at the other end, heuristic plausibility arguments, 
with \textbf{absolute certainty} as an unattainable limit point.

I've been publishing papers defending this thesis for more than a quarter of 
a century,\footnote
{See, for example, the introductory remarks in my 1974 \emph{J. ACM} paper [3].}
but few are convinced by my arguments.  So
in a recent paper [1] I've tried a new tactic. 
I use quotes from Leibniz, Einstein and G\"odel to make my case, like a lawyer
citing precedents in court\ldots
\vspace{1cm}

In spite of the fact that I regard the Riemann hypothesis as an excellent
new-axiom candidate---whether G\"odel agrees or merely thinks that a new axiom might be needed
to prove the RH, I'm not sure---let me briefly wax enthusiastic over a possible approach to a proof
of the RH that involves \textbf{randomness}. \emph{Disclaimer}: I'm not an expert on the RH.
What I'm about to relate is definitely an outsider's first impression, 
not an expert opinion.

\section*{A possible attack on the Riemann hypothesis?}

Here is a
concrete approach to the RH, one that uses no complex numbers. 
It's a 
probabilistic approach, and it involves the notion of \textbf{randomness}.
It's originally due to Stieltjes, who erroneously claimed to have proved the RH with
a variant of this approach.

The M\"obius $\mu$ function is about as likely to be 
$+1$ or $-1$ (see Derbyshire, \emph{Prime Obsession,} pp.\ 322--323).
\[
   \mu(n) =
\left\{
\begin{array}{ll}
   0 & \mbox{if $k^2$ divides $n$, $k > 1$,} 
\\
   (-1)^{\mbox{\scriptsize number of different prime divisors of $n$}}
   &
   \mbox{if $n$ is square-free.}
\end{array}
\right.
\]
The RH is equivalent to assertion that as $k$ goes from 1 to $n$, $\mu(k)$ 
is positive as often as negative.
More precisely, the RH is closely related to the assertion that 
the difference between 
\begin{itemize}
\item
the number of $k$ from 1 to $n$
for which $\mu(k) = -1$, and
\item
the number of $k$ from 1 to $n$
for which $\mu(k) = +1$ 
\end{itemize}
is $O(\sqrt n)$, of the order of square root of $n$, i.e., is bounded by a constant times
the square root of $n$.
This is roughly the kind of behavior that one would expect if the sign of the $\mu$
function \textbf{were chosen at random}
using independent tosses of a fair coin.\footnote
{For a more precise idea of what
to expect if the sign of the $\mu$
function were chosen at random, see the chapter on
the law of the iterated logarithm in Feller, 
\emph{An Introduction to Probability Theory and Its
Applications,} vol.\ 1, VIII.5 through VIII.7.}

This is usually formulated in terms of the Mertens function $M(n)$:\footnote{See [4, 5].}
\[
   M(n) = \sum_{k=1}^n \mu(k).
\]
According to Derbyshire, pp.\ 249--251, 
\[
   M(n) = O(\sqrt n) 
\]
implies the RH, but is actually stronger than the RH.  The RH is equivalent to the
assertion that for any $\epsilon > 0$, 
\[
   M(n) = O(n^{\frac{1}{2} + \epsilon}).
\]
\textbf{Could this formula be the door to the RH?!}

This probabilistic approach caught my eye while I was reading
this May's crop of RH books.  

I have always had an interest in probabilistic methods
in elementary number theory. This was one of the things that inspired me to come up with
my definition of \textbf{\emph{algorithmic} randomness} and to find 
algorithmic randomness in arithmetic [6]
in connection with diophantine equations.
However, I doubt that this work on algorithmic randomness is directly applicable to the RH.

In particular, these two publications greatly interested me as a child:
\begin{itemize}
\item
Mark Kac, 
\emph{Statistical Independence in Probability, Analysis and Number Theory,}
Carus Mathematical Monographs, vol.\ 12, Mathematical Association of America, 1959. 
\item
George P\'olya, ``Heuristic reasoning in the theory of numbers,'' 1959, reprinted in
Gerald W. Alexanderson, 
\emph{The Random Walks of George P\'olya,}
Mathematical Association of America, 2000. 
\end{itemize}
I think that anyone contemplating a probabilistic attack on the RH via the $\mu$ function
should read these two publications.
There is also some interesting work on random sieves, 
which are probabilistic versions of the sieve of Eratosthenes:
\begin{itemize}
\item
D. Hawkins, ``Mathematical sieves,'' \emph{Scientific American,} December 1958, pp.\ 105--112.
\end{itemize}

As P\'olya shows in the above paper---originally 
\emph{American Mathematical Monthly} \textbf{66}, pp.\ 375--384---probabilistic 
heuristic reasoning can do rather well
with the distribution of twin primes.  
By the way, this involves
Euler's $\gamma$ constant.
\textbf{Can a refinement of P\'olya's technique shed new light on $\mu$ and
on the RH?} I don't know, but I think that this is an interesting possibility.

By the way,
$\mathbf{P} \neq \mathbf{NP}$
also involves randomness, for as 
Charles Bennett
and John Gill showed in 
1981---\emph{SIAM Journal on Computing} \textbf{10}, pp.\ 96--113---with
respect (relative) to a random oracle $A$, 
$\mathbf{P}^A \neq \mathbf{NP}^A$
with probability one [7].

\section*{Further reading---Four ``subversive'' books}

\begin{itemize}
\item
On experimental mathematics:
\\  
Borwein, Bailey and Girgensohn, 
\emph{Mathematics by Experiment,} 
\emph{Experimentation in Mathematics,} A. K. Peters, in press.
\\   
(See [8]. There is a chapter on zeta functions in volume two.)
\item
On a quasi-empirical view of mathematics:
\\   
Tymoczko, 
\emph{New Directions in the Philosophy of Mathematics,} Princeton University Press, 1998.
\item
On pragmatically-justified new axioms and information-theoretic incompleteness:
\\   
Chaitin, 
\emph{From Philosophy to Program Size,} Tallinn Cybernetics Institute, 2003.
\\
(There is also an electronic version of this book [2].)
\end{itemize}

And regarding the adverse reaction of the mathematics community to the
ideas in the above four books, I think that it is interesting to recall
G\"odel's difficulties at the Princeton Institute for Advanced Study, as recounted in:
\begin{itemize}
\item
John L. Casti, 
\emph{The One True Platonic Heaven,}
John Henry Press, 2003.
\end{itemize}
According to Casti, 
one of the reasons that it took so long for G\"odel's appointment at the IAS to be converted
from temporary to
permanent is that
some of G\"odel's colleagues dismissed his incompleteness theorem.  
Now of course G\"odel has become a cultural icon\footnote
{In this connection, I should mention \emph{Incompleteness, a play and a theorem} by Apostolos 
Doxiadis, which is a play about G\"odel.  For more information, see [9].}
and mathematicians take 
incompleteness more seriously---but perhaps \textbf{not seriously enough}.

Mathematicians shouldn't be cautious lawyers---I much prefer the bold Eulerian way of doing
mathematics.
Instead of endlessing polishing, how about some adventurous pioneer spirit?
Truth can be reached through successive approximations; insistence on instant absolute rigor is 
sterile---that's what I've learned from incompleteness.\footnote
{In this connection, see da Costa, French, \emph{Science and Partial Truth.}}

\section*{Web References}

{\footnotesize
\begin{itemize}
\item[{[1]}] 
Two philosophical applications of algorithmic information theory.\\
\texttt{http://www.cs.auckland.ac.nz/CDMTCS/chaitin/dijon.html}
\item[{[2]}] 
From philosophy to program size.\\
\texttt{http://www.cs.auckland.ac.nz/CDMTCS/chaitin/ewscs.html}
\item[{[3]}] 
Information-theoretic limitations of formal systems.\\
\texttt{http://www.cs.auckland.ac.nz/CDMTCS/chaitin/acm74.pdf}
\item[{[4]}] 
Mertens function.\\
\texttt{http://mathworld.wolfram.com/MertensFunction.html}
\item[{[5]}] 
Mertens conjecture.\\
\texttt{http://mathworld.wolfram.com/MertensConjecture.html}
\item[{[6]}] 
Randomness in arithmetic.\\
\texttt{http://www.cs.auckland.ac.nz/CDMTCS/chaitin/sciamer2.html}
\item[{[7]}] 
Relative to a random oracle $A$, $\mathbf{P}^A \neq \mathbf{NP}^A \neq \mbox{co-}\mathbf{NP}^A$ 
with probability 1.\\
\texttt{http://www.research.ibm.com/people/b/bennetc/bennettc1981497f3f4a.pdf}
\item[{[8]}] 
Experimental mathematics website.\\
\texttt{http://www.expmath.info}
\item[{[9]}] 
Apostolos Doxiadis home page.\\
\texttt{http://www.apostolosdoxiadis.com}
\end{itemize}
}

\end{document}